\documentclass{amsart}

\usepackage{amssymb}
\usepackage{amscd}
\usepackage{amsmath}
\usepackage{amsfonts}

\newtheorem{theorem}{Theorem}
\newtheorem{prop}{Proposition}
\newtheorem{lemma}{Lemma}

\newtheorem*{ac}{Artin's Conjecture}
\newtheorem*{sac}{Strong Artin Conjecture}
\numberwithin{equation}{section} \numberwithin{lemma}{section}
\numberwithin{prop}{section}

\newcommand{\im}{\mathrm{Im}}
\newcommand{\C}{{\mathbb C}}

\newcommand{\Q}{{\mathbb Q}}
\newcommand{\A}{{\mathbb A}}

\newcommand{\Fr}{\mathrm {Fr}}
\newcommand{\G}{\mathrm {Gal}}
\newcommand{\GL}{\mathrm {GL}}
\newcommand{\PGL}{\mathrm {PGL}}
\newcommand{\PGSp}{\mathrm {PGSp_4}}

\newcommand{\SO}{\mathrm {SO}}
\newcommand{\GSp}{\mathrm {GSp_4(\C)}}
\newcommand{\GO}{\mathrm {GO_4(\C)}}
\newcommand{\Sp}{\mathrm {Sp_4}}
\newcommand{\PSp}{\mathrm {PSp_4}}

\newcommand{\s}{\sigma}

\newcommand{\zf}{\zeta_5}
\newcommand{\La}{\Lambda}
\newcommand{\ext}{\La^2}
\newcommand{\iso}{\simeq}

\newcommand{\sdp}{\rtimes}

\newcommand{\corr}{\leftrightarrow}

\newcommand{\gpi}{E_{2^4} \sdp C_5}

\title{A symplectic case of Artin's Conjecture}

\begin{document}

\maketitle \vspace{-0.2 in}
\begin{center}
{\sc Kimball Martin} \\ 253-37 Caltech \\ Pasadena, CA 91125 \\
{\em kimball@caltech.edu} \\ June 13, 2003 \vspace{0.1in}
\end{center}

\section{Introduction}

Here we will prove a new 4-dimensional symplectic case of Artin's
conjecture.  Let us first recall Artin's conjecture.  Let $G$ be
the Galois group of a finite Galois extension of number fields.
Let $\rho$ be a representation of $G$ over $\C$. In this context
one can associate to $\rho$ a meromorphic function $L(s,\rho)$
called the Artin $L$-function associated to $\rho$.

\begin{ac}
If $\rho$ does not contain the trivial representation, then
$L(s,\rho)$ is entire.
\end{ac}

An equivalent form of this conjecture is the following.  If $\rho$
is irreducible and non-trivial, then $L(s,\rho)$ is entire.  Artin
proved this conjecture in the case where $\rho$ is 1-dimensional
using his reciprocity law together with a result of Hecke. As
$L$-functions are inductive, this proved the conjecture if $\rho$
is {\em monomial}, i.e. induced from a degree one representation
of some subgroup.

Now consider the case where $\rho$ is 2-dimensional.  Let
$\bar{\rho}: G \to \PGL_2(\C)$ be the composition of $\rho: G \to
\GL_2(\C)$ with the natural projection of $\GL_2(\C)$ to
$\PGL_2(\C)$.  Let $\bar{G}$ denote the image of $G$ under
$\bar{\rho}$.  So $\bar{G}$ is a finite subgroup of $\PGL_2(\C)$,
which is isomorphic to $\SO_3(\C)$.  The only finite subgroups of
$\SO_3(\C)$ are cyclic, dihedral, tetrahedral (isomorphic to
$A_4$), octahedral (isomorphic to $S_4$) and icosahedral
(isomorphic to $A_5$). We thus classify $\rho$ according the
isomorphism type of $\bar{G}$. If $\bar{G}$ is cyclic or dihedral,
then $\rho$ is reducible or monomial and $L(s,\rho)$ is entire by
Artin's result.  Much later, Langlands applied to this problem the
theory of {\em automorphic representations}, which also have
associated $L$-functions.  For {\em cuspidal} automorphic
representations of $\GL_n$ the associated $L$-function {\em is} known to be
entire \cite{Ja}. Both automorphic and Artin $L$-functions can be
written as Euler products of {\em local factors} $L(s,\pi) =
\prod_v L(s,\pi_v)$ and $L(s,\rho) = \prod_v L(s,\rho_v)$.
Langlands formulated the following amazing conjecture
\cite{Lalet}.

\begin{sac}
Let $G$ be the Galois group of an extension $K/F$ of number
fields.  Let $\rho$ be an $n$-dimensional complex representation
of $G$. There exists an automorphic representation $\pi$ of
$\GL_n(\A_F)$, such that the $L$-functions agree almost
everywhere, i.e. except at a finite number of places $v$,
$L(s,\rho_v) = L(s,\pi_v)$. Moreover, if $\rho$ is irreducible,
then $\pi$ is cuspidal.
\end{sac}

We say $\rho$ is {\em modular} if the strong Artin conjecture
holds for $\rho$. Furthermore, if $\pi$ is as in the conjecture,
we will say $\rho$ {\em corresponds to} $\pi$ and write $\rho
\corr \pi$ or $\pi \corr \rho$.

The strong Artin (also called Langlands' reciprocity) conjecture
is really a stronger statement than Artin's conjecture.  For
example, nilpotent and supersolvable groups are monomial and hence
their representations $\rho$ satisfy Artin's conjecture. Arthur
and Clozel recently proved (\cite{AC}) that if $\rho$ has
nilpotent image, $\rho$ is also modular. However it has not yet
been shown that all $\rho$ with supersolvable image must be
modular.

Langlands proved the strong Artin conjecture, and hence Artin's
conjecture, in the tetrahedral case \cite{La}, i.e. when $\bar{G}
\iso A_4$ (and $\dim \rho = 2$). In his proof, he used three
important elements: the symmetric square lift from $\GL_2$ to
$\GL_3$ of Gelbart and Jacquet \cite{GJ}, {\em normal cubic base
change} on $\GL_2$ developed by him in \cite{La}, and the
structure of the group $A_4$. After non-normal cubic base change
was proven by Jacquet, Piatetski-Shapiro and Shalika \cite{JPSS},
Tunnell used it to extend Langlands' argument and prove the strong
Artin conjecture for the octahedral case ($\bar{G} \iso S_4$)
\cite{Tu}. This completed the solvable case in dimension 2. The
non-solvable (icosahedral) case is much more difficult. However,
partial but outstanding progress has been made recently in this
case by Buzzard, Dickinson, Shepherd-Barron and Taylor
\cite{BDST}.  For modular 2-dimensional $\rho$, the corresponding
symmetric $m$-th power representations for $m \leq 4$ are also
modular by the work of Kim and Shahidi \cite{KS}, \cite{Ki}.

There are not many other cases of Artin's conjecture known in
higher dimensions.  Progress has been primarily made only for
essentially self-dual (i.e., orthogonal or symplectic)
representations. Let $\rho$ be an essentially self-dual
irreducible complex $n$-dimensional representation of $G$ with
solvable image. If $n$ is odd, then $\rho$ maps into
$\mathrm{GO}_n(\C)$ and $\rho$ is monomial \cite{Ra}.  Therefore
in the odd-dimensional case Artin's conjecture (but not strong
Artin) is known.  If $n$ is even, $\rho$ has image either in
$\mathrm {GO}_n(\C)$ or $\mathrm {GSp}_n(\C)$. In the case $\rho:
G \to \GO$ has solvable image, then Ramakrishnan recently showed
$\rho$ satisfies the strong Artin conjecture by proving the
modularity of the Asai representation \cite{Ra}.

Here we make progress in the case $\rho$ has image in $\GSp$. We
are able to exhibit another group whose 4-dimensional
representations can be proven modular in a manner similar to
Langlands' tetrahedral argument. Consider $\rho: G \to \GSp$. Then
$\bar{G} \subseteq \PGSp(\C) = \GSp/\C^*$, which is isomorphic to
$\SO_5(\C)$.  The finite subgroups of $\PGSp(\C)$ have been
classified in \cite{Mi} and \cite{CM}.  Let $E_{2^4}$ denote the
elementary abelian group of order $2^4$ and $C_5$ the cyclic group
of order 5.  Then there is a semidirect product $\gpi$ with $C_5$
acting fixed point freely on $E_{2^4}$ contained inside
$\PGSp(\C)$.  If $\bar{G} \iso \gpi$, we show $\rho$ is modular.
Our proof uses the recent construction of Kim \cite{Ki} of the
exterior square $\ext(\pi)$ of an (isobaric) automorphic
representation $\pi$ of $\GL_4$, which agrees locally almost
everywhere with an automorphic representation of $\GL_6$.  We also
crucially use {\em normal quintic base change} (\cite{AC}) and the
structure of our group $\gpi$. An investigation of the other
subgroups of $\PGSp(\C)$ will appear in the author's thesis.  We now state the main result precisely.

\begin{theorem}
Let $L/F$ be a Galois extension of number fields and $\rho$ be an
irreducible 4-dimensional representation of $G = \G(L/F)$ into
$\GSp$. Suppose $\bar{G} = {\mathrm Im}(\bar{\rho}) \iso \gpi$.
Then $\rho$ is modular.
\end{theorem}

One pre-image of $E_{2^4} \sdp C_5$ inside $\Sp(\C)$ is $G =
E_{2^5} \sdp C_5$, where $E_{2^5}$ is the extraspecial group of
order $2^5$ isomorphic to the central product of $Q_8$ and $D_8$ with
identified centers (see \cite{Do} for definitions).  Consider an irreducible 4-dimensional
complex representation $\rho$ of $G$.  As $G$ has no subgroups of
index 2 or 4, $\rho$ is primitive, i.e. not induced.  Thus Artin's
conjecture for $\rho$ does not follow from previous results.

We know that examples of $E_{2^5} \sdp C_5$ extensions of $\Q$ exist by Shafarevich's Theorem (\cite{NSW})  because $E_{2^5} \sdp C_5$ is solvable.
Though Shafarevich's proof is non-constructive, we can illustrate how to construct such an extension in our case.  Let $\alpha_i = \zeta_{11}^i + \zeta_{11}^{-i}$, where $\zeta_{11}$ is a complex $11$-th root of unity.  Let $E$ be the cyclic Galois extension $\Q(\alpha_1)$ of $\Q$
of degree 5.  It is known how to construct $Q_8$ and $D_8$ extensions of a number field  (\cite{JLY}).  Let $K = E(\sqrt{1+A+B+AB})$ and $M=E(\sqrt{\alpha_1+i\alpha_2},\sqrt{\alpha_2+\alpha_4+4}),$ where $A=(3+\alpha_5)^{-1/2}$
and $B=(1+\alpha_1^2+\alpha_1^2\alpha_3^2)^{-1/2}$.  Then $\G(K/E) \iso Q_8$ and $\G(M/E) \iso D_8$.  The compositum $KM$ has three normal (over $E$) subextensions of index 2.  Let $L/E$ be the one corresponding to the central product of $Q_8$ with $D_8$.  Then $L/\Q$ is Galois with
Galois group $E_{2^5} \sdp C_5$.

Before we go on, let us be a little more explicit about the
similarities between Langlands' tetrahedral case and Theorem 1. We
will first outline Langlands' argument.  Let $\rho$ be a
representation of $G$ into $\GL_2(\C)$ with image $A_4$.  Then $G$
has a normal subgroup $H$ of index 3.  Then $H$ is dihedral.  All
representations of dihedral groups are modular (\cite{La}), so
$\rho_H$ is modular. Let $\Pi$ be an automorphic representation
corresponding to $\rho_H$. Normal cubic base change for $\GL_2$
tell us that there are three representations $\pi_i$ whose {\em
base change} $\pi_{i,H}$ (the automorphic version of restricting a
representation) to $H$ is $\Pi$.  So we think one of these $\pi_i$
should correspond to $\rho$.  Using the determinant of $\rho$,
Langlands proved there is a unique $\pi$ among the $\pi_i$'s with
its central character $\omega_\pi$ corresponding to $\det(\rho)$,
such that $Sym^2(\rho) \corr Sym^2(\pi)$. Langlands then showed
that if $\rho_H \corr \pi_H$ and $Sym^2(\rho) \corr Sym^2(\pi)$,
but $\rho$ does not correspond to $\pi$, then $\bar{G} \iso A_4$
must have an element of order 6. But $A_4$ has no elements of
order 6, therefore $\rho$ must correspond to $\pi$, i.e. $\rho$ is
modular.

The basic approach to proving Theorem 1 is similar to the
tetrahedral argument, but with several differences.  The first
difficulty encountered in our case is that the determinant alone
does not give enough information to choose the appropriate $\pi$.
More precisely, suppose $\rho$ is a representation of $G$ into
$\GSp$ with projective image $\gpi$.  Then $G$ has a normal
subgroup $H$ of index 5. Since $H$ is nilpotent, $\rho_H$ is
modular. Say $\rho_H \corr \Pi$.  Normal quintic base change tells
us that there are five representations $\pi_i$ whose base change
$\pi_{i,H}$ to $H$ equals $\Pi$.  Using the fact that $\ext(\rho)$
has an invariant line, we are able to pick a unique $\pi$ from the
$\pi_i$'s such that $\ext(\rho) \corr \ext(\pi)$. Then we deduce
either $\rho \corr \pi$ or $\bar{G}$ has an element of order 10,
which it does not.  Therefore $\rho$ must be modular.

In fact, Langlands' conjectures actually predict $\rho$ is {\em modular of symplectic type},
i.e $\rho$ corresponds to an automorphic representation of $\mathrm{GSp}_4(\A_F)$.  We discuss this briefly in the final section.

\section{Preliminaries}
\numberwithin{theorem}{section}

In this section, we will review the theories of $L$-functions and
base change.  Those already familiar with them may wish to skip to
the proof in the next section.

First we will begin with some notation.  Let $F$ be a number field
and $L$ be a finite Galois extension of $F$ with Galois group $G$.
Let $\rho$ be a representation of $G$ into $\GL_n(\C)$.  We will
refer to restriction and induction of representations with the
corresponding fields. More precisely, let $E$ be a subextension of
$L/F$ and $H = \G(L/E)$.  Then $\rho_E$ denotes the restriction
$\rho_H$ of $\rho$ to the subgroup $H$.  For a representation $\s$
of $H$, $I_E^F \s$ denotes the induced representation $I_H^G \s$
of $\s$ from $H$ to $G$.  We will use $\frak N_{E/F}$ to denote
the norm map from $E$ to $F$.

Let $v$ be a place of $F$ and $\Fr_v$ be the corresponding
Frobenius class.  Denote by $q_v$ the size of the residue field
${\mathcal O}_F/v$.  Recall the Artin $L$-function is defined on
some right-half plane by $L(s,\rho) = \prod_v L(s,\rho_v)$ where
the product is over all places $v$ of $F$ and if $v$ is
unramified,
\[ L(s,\rho_v) = \frac{1}{det(1- \rho(\Fr_v) q_v^{-s})}. \]
Brauer showed for any $\rho$, $L(s,\rho)$ extends to a meromorphic
function on all of $\C$.  In the case $\rho$ is the trivial
representation {\pmb 1}, the $L$-function is the Dedekind zeta
function $\zeta_F(s)$ of $F$.  For a complete definition of Artin
$L$-functions, see \cite{Ro}, \cite{Ne} or \cite{Ma}.

One can also define $L$-functions for automorphic representations.
Let $\A_F$ be the adeles of the number field $F$.

\begin{theorem}\label{lthm} {\rm (\cite{Ja})}
Let $\pi$ be an automorphic representation of $\GL_n(\A_F)$. Then
we can associate local factors $L(s,\pi_v)$ to each $v$ and define
a holomorphic $L$-function $L(s, \pi) = \prod L(s,\pi_v)$ in some
right-half plane.  Moreover $L(s,\pi)$ extends to an meromorphic
function on $\C$, which is actually entire if $\pi$ is cuspidal
and non-trivial.
\end{theorem}

For more information on automorphic representations and their
$L$-functions, see \cite{Geintro}, \cite{Kn}, or \cite{Gebook}.
The strong Artin conjecture asserts that given a Galois
representation $\rho$ as before, there exists an automorphic
representation $\pi$ corresponding to $\rho$ (see Introduction).
If $\rho \corr \pi$ and $\rho$ is irreducible, then the following
result implies $\pi$ is cuspidal.

\begin{theorem} {\rm (\cite{JS})} \label{js}
Let $\pi$ be an automorphic representation of $\GL_n(\A_F)$ and
$\check \pi$ its {\em contragredient}. Then $L(s, \pi \times
\check \pi)$ has a simple pole at $s=1$ if and only if $\pi$ is
cuspidal.
\end{theorem}

For if $\rho$ is irreducible and $\check \rho$ its contragredient,
then $L(s, \rho \otimes \check \rho)$ has a simple pole at $s=1$.
Moreover if $\rho \corr \pi$, then $L(s, \rho \otimes \check
\rho)$ and $L(s, \pi \times \check \pi)$ agree almost everywhere
locally.  So by Theorem \ref{js}, $\pi$ is indeed cuspidal.

Now to show that the strong Artin conjecture actually implies
Artin's conjecture we only need the following standard fact.

\begin{prop}
If $\pi$ is cuspidal and $L(s, \pi_v) = L(s,\rho_v)$ for almost
all $v$, then in fact $L(s, \pi) = L(s, \rho)$.
\end{prop}

So if $\rho$ is irreducible and non-trivial and corresponds to
$\rho$, then by Theorem \ref{js} $\rho$ is cuspidal (and
non-trivial). The proposition tells us that in fact $L(s, \pi) =
L(s, \rho)$.  By Theorem \ref{lthm} we know that $L(s,\pi)$ is
entire, i.e. $L(s, \rho)$ is entire, i.e. Artin's conjecture is
true for $\rho$.

If $\pi$ and $\pi'$ are cuspidal representations of $\GL_m(\A_F)$
and $\GL_n(\A_F)$ respectively, then one can form their {\em
isobaric sum} $\pi \boxplus \pi'$ which is an automorphic
representation of $\GL_{m+n}(\A_F)$ \cite{JS}.  If $\rho \corr
\pi$ and $\rho' \corr \pi'$, then $\rho \oplus \rho' \corr \pi
\boxplus \pi'$.  So every automorphic representation which
corresponds to a Galois representation will be {\em isobaric},
i.e. a finite isobaric sum of cuspidal representations.

A powerful tool to prove the strong Artin conjecture in certain
instances is the theory of base change. Base change is an
operation on automorphic representations which corresponds to the
restriction of Galois representations.  We now list the important
properties of base change that we need.

\begin{theorem}\label{bcthm} {\rm (\cite{AC})}
Let $L/F$ be a Galois extension of number fields.  Let $E/F$ be a
normal cyclic subextension of prime degree.  For each isobaric
representation $\pi$ of $\GL_n(\A_F)$, there exists a unique automorphic
representation of $\GL_n(\A_E)$ called the base change of $\pi$ to
$E$ and denoted by $\pi_E$ such that:

(i) (descent) a cuspidal representation $\Pi$ of $\GL_n(\A_E)$ is
the base change $\pi_E$ of some $\pi$ if and only if $\Pi$ is {\em
Galois invariant} (in particular, if $\Pi \corr \rho_E$ where
$\rho$ is some representation of $\G(L/F)$);

(ii) if $\pi'$ is also an isobaric representation of $\GL_n(\A_F)$
then $\pi_E = \pi'_E$ if and only if $\pi' = \pi \otimes \delta$
for some idele class character $\delta$ of $F^* \frak
N_{E/F}(\A^*_E) \backslash \A^*_F \iso \G(E/F)$;

(iii) (compatibility with reciprocity) if $\rho$ is a
representation of $\G(L/F)$, $\rho \corr \pi$ and $\rho_E$ is
modular, then $\rho_E \corr \pi_E$; and

(iv) (compatibility with twisting) if $\chi$ is an idele class
character of $F$ and $\chi_E = \chi \circ N_{E/F}$, then
\[ ( \pi \otimes \chi)_E = \pi_E \otimes \chi_E. \]
\end{theorem}

The complementary construction to base change is {\em automorphic
induction}, which corresponds to induction of Galois
representations.

\begin{theorem} \label{ai} {\rm (\cite{AC},\cite{HH})}
Let $L/F$ be a Galois extension of number fields.  Let $E/F$ be a
normal cyclic subextension of prime degree.  Let $\rho$ be a
complex representation of $\G(L/E)$ and suppose $\rho \corr \pi$,
for some automorphic representation $\pi$. Then there exists an
{\em induced} automorphic representation, denoted $I_E^F\pi$ such
that $I_E^F \rho \corr I_E^F\pi$.
\end{theorem}

\section{Proof of Theorem 1}

Let $L/F$ be a (finite) Galois extension of number fields with
Galois group $G$.  Suppose $\rho$ is an (injective) representation of $G$ into
$\GSp$ such that $\bar{G} \iso \gpi$.  Let $E$ be the normal
quintic subextension of $L/F$ corresponding to the pre-image of
$E_{2^4}$.

\begin{lemma} The representations $\rho_E$ and
$\ext(\rho)$ are modular.
\end{lemma}

\begin{proof}
As $\G(L/E)$ is a cyclic central extension of a 2-group, it is a
direct product of a 2-group $P_2$ with a cyclic group $C$ of odd
order. Therefore $\G(L/E)$ is nilpotent.  By a theorem of Arthur
and Clozel, all representations of nilpotent groups are modular
\cite{AC}. In particular $\rho_E$ is modular.

Since $\rho$ is of symplectic type, $\ext(\rho)$ has an invariant
line.  Write $\ext(\rho) = \nu \oplus r$ where $\nu$ is
1-dimensional and $r$ is  5-dimensional.  Note $r$ is irreducible
because it factors through $E_{2^4} \sdp C_5$, which only has 1-
and 5- dimensional irreducible representations.  We claim $r$ is
induced from $E$.  As $P_2$ is a 2-group, every irreducible
representation of $P_2$ has dimension $2^j$ for some $j$.
Therefore the same is true for $\G(L/E) \iso P_2 \times C$. Hence
in the decomposition of $r_E$ into its irreducible components, we
must have a 1-dimensional representation $\lambda$. In particular,
$r = I_E^F \lambda$. Since $E$ is a normal subextension $r$ is
modular by Theorem \ref{ai}, whence $\ext(\rho)$ is also.
\end{proof}

Let us say $\rho_E \corr \Pi$.  We claim $\rho_E$ is irreducible.
Indeed, $\rho$ irreducible implies that $\G(E/F) \iso C_5$ acts
transitively on the irreducible components of $\rho_E$.  This
action has order dividing 5.  Thus if there is more than one
irreducible component of $\rho_E$, there must be five or a
multiple thereof.  However $\dim \rho_E = 4$, so that is
impossible.  Then by Theorem \ref{js}, $\Pi$ is cuspidal.

We can apply Theorem \ref{bcthm}(i) to get an automorphic
representation $\pi_0$ of $\GL_4(\A_F)$ whose base change
$\pi_{0,E}$ corresponds to $\rho_E$.  Let $\delta = \delta_{E/F}$
be a non-trivial idele class character of $F^* \frak
N_{E/F}(\A^*_E) \backslash \A^*_F \iso \G(E/F) \iso C_5$.  Let
$\pi_i = \pi_0 \otimes \delta^i$ for $i = 1, 2, ..., 4$.  Since
$\delta_E = {\pmb 1}$, all the $\pi_i$'s base change to $\pi_{i,E}
\iso (\pi_0 \otimes \delta)_E \iso \pi_{0,E} \iso \Pi$ by Theorem
\ref{bcthm}(iv).  In fact, part (ii) of the same theorem tells us
that these are all the cuspidal representations of $\GL_4(\A_F)$
whose base change to $E$ is $\Pi$.

\begin{lemma} There is a unique $i \in \{ 0, 1, 2, ..., 4 \}$ such that
$\ext(\pi_i) \corr \ext(\rho)$.
\end{lemma}

\begin{proof}
All the representations $\ext(\pi_i)$ base change to $\ext(\pi_0
\otimes \delta^i)_E = \ext(\pi_0)_E$ by Theorem \ref{bcthm}(iv).
They are all distinct because they have distinct central characters
$\omega_{\ext(\pi_i)} = \omega_{\ext(\pi_0)} \delta^{2i}$.

Theorem \ref{bcthm}(ii) then yields that the $\ext(\pi_i)$ are the
only representations which base change to $\ext(\pi_0)_E$. But by
part (iii) of this theorem, the automorphic representation on
$\GL_6(\A_F)$ which corresponds to $\ext(\rho)$ must also base
change to $\ext(\pi_0)_E$.  Thus for some $i$, $\ext(\pi_i) \corr
\ext(\rho)$.
\end{proof}

Denote the $\pi_i$ of the lemma by $\pi$. We claim now that in
fact $\rho \corr \pi$. It will suffice to show for all
unramified places that $\rho_v \corr \pi_v$.  Say
$\rho_v$ has Frobenius eigenvalues $\{a,b,c,d\}$ and $\pi_v$ has Satake
parameters
$\{e,f,g,h\}$.  We want to show
$\{a,b,c,d\} = \{e,f,g,h\}$.
Let $D$ be a diagonal element of $\GL_4$. Then
$\ext(D) = 1$ if and only if $D = \pm I$. Hence $\ext(\rho_v)
\corr \ext(\pi_v)$ implies
\begin{equation}
\{a,b,c,d\} = \pm \{e,f,g,h\}.
\end{equation}

If they are equal, we are done.  Assume therefore
\begin{equation}
\{a,b,c,d\} = -\{e,f,g,h\}.
\end{equation}
Now we can use base change to $E$.  In our group $\bar{G}$ any
element raised to the 5th power lies inside $E_{2^4}$ (see Lemma
\ref{gp} below). Thus any element of $G$ raised to the 5th power
lies inside $\G(L/E)$, the pre-image of $E_{2^4}$.  In particular
$Fr_v^5 \in {\mathcal O}_{E_w}$, where $w$ is a prime of $E$ above
$v$. Then $\rho_{v,E} \corr \pi_{v,E}$ implies
$\{a^5,b^5,c^5,d^5\} = \{e^5,f^5,g^5,h^5\}$.  By our assumption we
have,
\begin{equation}
\{a^5,b^5,c^5,d^5\} = \{-a^5,-b^5,-c^5,-d^5\}.
\end{equation}
Without loss of generality, assume $a^5 = -b^5$ and $c^5 = -d^5$.
Then either $b = -\zf a$ or $d = -\zf c$, for otherwise $a = -b, c
= -d$ which would imply $\{a,b,c,d\} = \{e,f,g,h\}$.  Let us say
\begin{equation}
b = -\zf a.
\end{equation}

Then
\begin{equation}
 \rho(Fr_v) \sim
\left(\begin{array}{cccc}
a & 0 & 0 & 0 \\
0 & -\zf a & 0 & 0 \\
0 & 0 & c & 0 \\
0 & 0 & 0 & d \\
\end{array} \right),
\end{equation}
so \begin{equation}
 \bar{\rho}(Fr_v) \sim  \left(\begin{array}{cccc}
1 & 0 & 0 & 0 \\
0 & -\zf & 0 & 0 \\
0 & 0 & c/a & 0 \\
0 & 0 & 0 & d/a \\
\end{array} \right),
\end{equation}
is an element of order divisible by 10 in $\bar{G} =
Im(\bar{\rho})\subseteq \PSp(\C)$.  But $\bar{G}$ has no elements
of order 10 by Lemma \ref{gp} below, so $\rho$ is modular by
contradiction.

\begin{lemma} \label{gp}
Every element $g \in \gpi$ has order 5 except for the elements in
the normal subgroup $E_{2^4}$, which have order 1 or 2.
\end{lemma}

\begin{proof}
Let $g \in \gpi$ such that $g \not \in E_{2^4}$.  We can write $g
= az$ where $a \in E_{2^4}$ and $z \in C_5$, $z \neq 1$.  We claim
$g^5$ commutes with $z$.  Write
\begin{eqnarray}
g^5 &=& (az)(az)(az)(az)(az) \\
&=& a(zaz^{-1})(z^2az^{-2})(z^3az^{-3})(z^4az^{-4}). \label{g5}
\end{eqnarray}
We also have
\begin{eqnarray}
zg^5z^{-1} &=& (za)(za)(za)(za)(za) \\
&=& (zaz^{-1})(z^2az^{-2})(z^3az^{-3})(z^4az^{-4})a. \label{zg5z}
\end{eqnarray}
Each $z^jaz^{-j}$ lies in the normal abelian subgroup $E_{2^4}$ and
therefore the $z^jaz^{-j}$'s commute.  Thus we can rearrange the
terms in (\ref{zg5z}) to get (\ref{g5}) and we have
$zg^5z^{-1}=g^5$.

Now, since each term in (\ref{g5}) lies in $E_{2^4}$, then $g^5
\in E_{2^4}$ also.  But the action of $C_5$ on $E_{2^4}$ fixes
only the identity.  Thus $g^5 = 1$.
\end{proof}

\section{Transfer to $\mathrm{GSp}_4$}
As before, consider a Galois group $G$ of an extension of number
fields $L/F$ and an irreducible 4-dimensional representation
$\rho$ of $G$ into $\GSp$.  Suppose that $\rho$ is modular, i.e.
$\rho$ corresponds to some cuspidal representation $\pi$ of
$\GL_4(\A_F)$. The fact that $\im(\rho) \subseteq \GSp$ implies
that $L(s, \ext(\rho) \otimes \nu^{-1})$ has a simple pole at
$s=1$ for a suitable 1-dimensional representation $\nu$ of $G$
(the ``polarization'').

This implies that the corresponding automorphic $L$-function
$L(s,\pi ; \ext \otimes \nu^{-1})$ admits a pole at $s=1$.  An
unpublished theorem of Jacquet, Piatetski-Shapiro and Shalika says
that, because of this pole, $\pi$ {\em transfers} to a generic irreducible cuspidal
automorphic representation $\Pi$ of $\mathrm{GSp}_4(\A_F)$ with
central character $\nu$ such that,
\[ L^S(s,\Pi) = L^S(s,\pi), \]
for any finite set of primes $S$ outside of which $\pi$ is unramified.  Here the $L$-function on
the left is the degree 4 $L$-function of $\Pi$ studied in
\cite{PS}; and if $L(s)= \prod_v L_v(s)$ is an Euler product, then
$L^S(s)$ denotes the incomplete $L$-function $\prod_{v\not \in S}
L_v(s)$. Thus $\rho$ in fact corresponds to the cuspidal
representation $\Pi$ of $\mathrm{GSp}_4(\A_F)$, i.e. $\rho$ is
{\em modular of symplectic type} as predicted by Langlands.

However we are not stressing this here because this theorem of Jacquet, Piatetski-Shapiro and Shalika remains unpublished.  We hope to go into more detail in the thesis.  The key point is that $\GL_4$ maps into the connected component of $\mathrm{GO}_6$, and $\pi$ gives rise to a cuspidal automorphic representation
$\pi'$ of $\mathrm{GO}_6(\A_F)^0$.  The desired $\Pi$ is obtained by the theta correspondence.  The obstruction to this transfer is the residue of the pole of $L(s,\pi; \ext \otimes \nu^{-1})$ at $s=1$ (see \cite{JS2}).  Finally the ongoing work of J. Arthur will give another proof, using the trace formula, of the existence of $\Pi$ and other members of its packet (see \cite{Ar} for his program).

\vspace{5pt}
\section*{acknowledgements}
The author would like to foremost thank his advisor, Dinakar
Ramakrishnan, for suggestions and guidance throughout this work.
He would also like to thank David Wales and Michael Aschbacher for
their useful words about symplectic groups.  The author is grateful to
Jim Cogdell, David Whitehouse and the referee for comments that
led to several minor improvements and corrections.  Finally, thanks go to
the GAP Group for their wonderful package as many group and
character computations were done in the initial stages of this
work using \cite{GAP}.

\vspace{5pt}
\section*{references}

\end{document}